\documentstyle{amsppt}
\tolerance 3000
\pagewidth{5.5in}
\vsize7.0in
\magnification=\magstep1
\widestnumber \key{AAAAAAAAAAAAA}
\topmatter
\author Alex Iosevich, Nets Katz, and Terry Tao
\endauthor
\thanks The research of Alex Iosevich is partially supported by the NSF Grant DMS00-87339.
The research of Nets Katz is partially supported by the NSF Grant
DMS-9801410. Terry Tao is a Clay prize fellow and supported by
grants from the Packard and Sloan foundations.
\endthanks
\title Fuglede conjecture holds for convex planar domains
\endtitle
\abstract Let $\Omega$ be a compact convex domain in the plane. We
prove that $L^2(\Omega)$ has an orthogonal basis of exponentials
if and only if $\Omega$ tiles the plane by translation.
\endabstract
\endtopmatter

\def\dist{\hbox{{\rm dist}}}

\def\R{\hbox{{\bf R}}}
\def\Z{\hbox{{\bf Z}}}
\def\eps{\varepsilon}

\document

\head Section 0: Introduction \endhead

Let $\Omega$ be a domain in ${\Bbb R}^d$, i.e., $\Omega$ is a
Lebesgue measurable subset of $\Bbb{R}^d$ with finite non-zero
Lebesgue measure. We say that a set $\Lambda \subset {\Bbb R}^d$
is a {\it spectrum} of $\Omega$ if ${\{e^{2\pi i x \cdot \lambda}
\}}_{\lambda \in \Lambda}$ is an orthogonal basis of
$L^2(\Omega)$.

\example{Fuglede Conjecture}(\cite{Fug74}) A domain $\Omega$
admits a spectrum if and only if it is possible to tile ${\Bbb
R}^d$ by a family of translates of $\Omega$.
\endexample

If a tiling set or a spectrum set is assumed to be a lattice, then
the Fuglede Conjecture follows easily by the Poisson summation
formula. In general, this conjecture is nowhere near resolution,
even in dimension one. However, there is some recent progress
under an additional assumption that $\Omega$ is convex. In
\cite{IKP99}, the authors prove that the ball does not admit a
spectrum in any dimension greater than one. In \cite{Kol99},
Kolountzakis proves that a non-symmetric convex body does not
admit a spectrum. In \cite{IKT00}, the authors prove that any
convex body in ${\Bbb R}^d$, $d>1$, with a smooth boundary, does
not admit a spectrum. In two dimensions, the same conclusion holds
if the boundary is piece-wise smooth and has at least one point of
non-vanishing curvature. The main result of this paper is the
following: \proclaim{Theorem 0.1} Let $\Omega$ be a convex compact
set in the plane. The Fuglede conjecture holds. More precisely,
$\Omega$ admits a spectrum if and only if $\Omega$ is either a
quadrilateral or a hexagon.
\endproclaim

Our task is simplified by the following result due to
Kolountzakis. See \cite{Kol99}. \proclaim{Theorem 0.2} Convex
non-symmetric subsets of ${\Bbb R}^d$ do not admit a spectrum.
\endproclaim

Thus, it suffices to prove Theorem 0.1 for symmetric sets. Recall
that a set $\Omega$ is symmetric with respect to the origin when
$x \in \partial \Omega$ if and only if $-x \in \partial \Omega$.

This paper is organized as follows. The first section deals with
basic properties of spectra. The second section is dedicated to
the properties of the Fourier transform of the characteristic
function of a convex set. In the third section we prove Theorem
0.1 for polygons, and in the fourth section we prove that any
convex set which is not a polygon does not admit a spectrum, thus
completing the proof of Theorem 0.1.

\head Section 1: Basic properties of spectra \endhead

Let
$$Z_{\Omega}= \left\{\xi \in \Bbb R^d: \hat
\chi_{\Omega}(\xi)=\int_\Omega e^{-2\pi i \xi \cdot x}\ dx=0
\right\}. \tag1.1$$ The orthogonality of a spectrum $\Lambda$
means precisely that
$$ \lambda - \lambda' \in Z_\Omega \hbox{ for all } \lambda, \lambda' \in
\Lambda, \lambda \neq \lambda'.  \tag1.2$$

It follows that the points of a spectrum $\Lambda$ are separated
in the sense that
$$ |\lambda-\lambda'| \gtrsim 1 \ \text{for all }
\lambda \not=\lambda', \ \ \lambda, \lambda' \in \Lambda.
\tag1.3$$ Here, and throughout the paper, $a \lesssim b$ means
that there exists a positive constant $C$ such that $a \leq Cb$.
We say that $a \approx b$ if $a \lesssim b$ and $a \gtrsim b$.

The following result is due to Landau. See \cite{Lan67}. Let
$$ D_R^{+}= \max_{x \in {\Bbb R}^n} \# \{ \Lambda \cap Q_R(x)\}, \tag1.4$$
where $Q_R(x)$ is a cube of sidelength $2R$ centered at $x$, and let
$$ D_R^{-}= \min_{x \in {\Bbb R}^n}\#\{ \Lambda \cap Q_R(x)\}. \tag1.5$$

Then
$$ \limsup_{R \rightarrow \infty} \frac{D_R^{\pm}}{{(2R)}^n} =
|\Omega|. \tag1.6$$

It is at times convenient to use the following related result. We
only state the special case we need for the proof of Theorem 0.1.
For a more general version see \cite{IosPed99}. \proclaim{Theorem
1.1} Let $\Omega$ be a convex domain in ${\Bbb R}^2$. Then there
exists a universal constant $C$ such that if
$$ R \ge C
\left(\frac{{|\partial \Omega|}}{|\Omega|}\right), \tag1.7$$ then
$$ \Lambda \cap Q_R(\mu) \not= \emptyset \tag1.8$$ for every
$\mu \in {\Bbb R}^2$, and any set $\Lambda$ such that
$E_{\Lambda}$ is an exponential basis for $L^2(\Omega)$, where
$Q_R(\mu)$ denotes the cube of sidelength $2R$ centered at $\mu$.
\endproclaim

The proof of Theorem 1.1, and the preceding result due to Landau,
are not difficult. Both proofs follow, with some work, from the
fact that $\Omega$ admits a spectrum $\Lambda$ if and only if
$$ \sum_{\Lambda} {|\widehat{\chi}_{\Omega}(x-\lambda)|}^2 \equiv 1,
\tag1.9$$ and some averaging arguments. To say that $\Omega$
admits a spectrum $\Lambda$ means that the Bessel formula
${||f||}_{L^2(\Omega)}^2=\sum_{\Lambda} {|\hat{f}(\lambda)|}^2$
holds. Since the exponentials are dense, it is enough to establish
such a formula with $f=e^{2 \pi i x \cdot\xi}$, which is precisely
the formula $(1.9)$.

\head Section 2: Basic properties of $\widehat{\chi}_{\Omega}$ and
related properties of convex sets
\endhead

Throughout this section, and the rest of the paper, $\Omega$
denotes a convex compact planar domain. The first two results in
this section are standard and can be found in many books on
harmonic analysis or convex geometry. \proclaim{Lemma 2.1}
$|\widehat{\chi}_{\Omega}(\xi)| \lesssim
\frac{diam{\Omega}}{|\xi|}$. Moreover, if $\Omega$ is contained in
a ball of radius $r$ centered at the origin, then $|\nabla
\widehat{\chi}_{\Omega}(\xi)| \lesssim \frac{r^2}{|\xi|}$.
\endproclaim

The lemma follows from the divergence theorem which reduces the
integral over $\Omega$ to the integral over $\partial \Omega$ with
a factor of $\frac{1}{|\xi|}$, and the fact that convexity implies
that the measure of the boundary $\partial \Omega$ is bounded by a
constant multiple of the diameter. The second assertion follows
similarly. \proclaim{Lemma 2.2} Suppose that $\xi$ makes an angle
of at least $\theta$ with every vector normal to the boundary of
$\Omega$. Then
$$ |\widehat{\chi}_{\Omega}(\xi)| \lesssim
\frac{1}{ \theta {|\xi|}^2}. \tag2.1$$

Moreover, if $\Omega$ is contained in a ball of radius $r$, then
$|\nabla \widehat{\chi}_{\Omega}(\xi)| \lesssim \frac{r}{| \theta
{|\xi|}^2|}$.
\endproclaim

To prove this, one can again reduce the integral to the boundary
while gaining a factor $\frac{1}{|\xi|}$. We may parameterize a
piece of the boundary in the form $\{(s, -\gamma(s)+c): a \leq s
\leq b\}$, where $\gamma$ is a convex function, and, without loss
of generality, $c=0$, $a=0$, $b=1$, and $\gamma(0)=\gamma'(0)=0$.
We are left to compute
$$ \int_{0}^{1} e^{i(s \xi_1-\gamma(s)
\xi_2)} J(s)ds, \tag2.2$$ where $J(s)$ is a nice bounded function
that arises in the application of the divergence theorem. The
gradient of the phase function $s \xi_1-\gamma(s) \xi_2$ is $\xi_2
\left( \frac{\xi_1}{\xi_2}-\gamma'(s) \right)$, and our assumption
that $\xi$ makes an angle of at least $\theta$ with every vector
normal to the boundary of $\Omega$ means that the absolute value
of this expression is bounded from below by $|\xi_2| \theta$.
Integrating by parts once we complete the proof in the case
$|\xi_1| \lesssim |\xi_2|$. If $|\xi_1|>>|\xi_2|$, the absolute
value of the derivative of $s \xi_1-\gamma(s) \xi_2$ is bounded
below by $|\xi_1|$, so integration by parts completes the proof.
The second assertion follows similarly.

\proclaim{Lemma 2.3} Let $f$ be a non-negative concave function on
an interval $[-1/2,1/2]$. Then, for every $0<\delta \lesssim 1$,
there exists $R \approx \frac{1}{\delta}$ such that $|\hat{f}(R)|
\gtrsim \delta f \left( \frac{1}{2}-\delta \right)$.
\endproclaim

To see this, let $\phi$ be a positive function such that $\phi(x)
\lesssim {(1+|x|)}^{-2}$, $\widehat{\phi}$ is compactly supported,
and $\phi(0)=1$ in a small neighborhood of the origin. Consider
$$ \int f \left( \frac{1}{2}-\delta t
\right)(\phi(t+1)-K\phi(K(t+1)))dt, \tag2.3$$ where $f$ is defined
to be $0$ outside of $[a,b]$ and $K$ is a large positive number.
If $K$ is sufficiently large, $(\phi(t+1)-K\phi(K(t+1)))$ is
positive for $t>0$, and $\approx 1$ on $[\frac{1}{2},1]$. It
follows that
$$ \int f \left( \frac{1}{2}-\delta t
\right)(\phi(t+1)-K\phi(K(t+1)))dt \gtrsim f \left(
\frac{1}{2}-\delta \right). \tag2.4$$

Taking Fourier transforms, we see that
$$ \int \frac{1}{\delta}
\hat{f} \left( \frac{r}{\delta} \right) e^{i \pi
r}\left(\hat{\phi}(r)-\hat{\phi} \left( \frac{r}{K} \right)
\right) dr \gtrsim  f \left( \frac{1}{2}-\delta \right). \tag2.5$$

Multiplying both sides by $\delta$ and using the compact support
of $\hat{\phi}(r)-\hat{\phi} \left( \frac{r}{K} \right)$, we
complete the proof.

\proclaim{Corollary 2.4} Let $\Omega$ be a convex body of the
form
$$ \Omega=\{(x,y): a \leq x \leq b, \ -g(x) \leq y \leq
f(x)\}, \tag2.6$$ where $f$ and $g$ are non-negative concave
functions on $[a,b]$. Then for every $0<\delta \lesssim b-a$,
there exists $R \approx \frac{1}{\delta}$ such that
$$ |\widehat{\chi}_{\Omega}| \gtrsim \delta \left(f \left(
\frac{1}{2}-\delta \right)+g \left( \frac{1}{2}-\delta \right)
\right). \tag2.6$$
\endproclaim

\head Section III: Lattice properties of spectra \endhead
\vskip.125in

Let $\Omega$ be a compact convex body in $\R^2$ which is symmetric
around the origin, but is not a quadrilateral.  Let $\Lambda$ be a
spectrum of $\Omega$ which contains the origin. The aim of this
section is to prove the following two propositions which show that
if a spectrum exists, it must be very lattice-like in the
following sense.

\proclaim{Proposition 3.1}  Let $I$ be a maximal closed interval
in $\partial \Omega$ with midpoint $x$.  Then
$$ \xi \cdot 2x \in \Z \tag3.1$$ for all $\xi \in \Lambda$.
\endproclaim

\proclaim{Proposition 3.2}  Let $x$ be an element of $\partial
\Omega$ which has a unit normal $n$ and which is not contained in
any closed interval in $\Omega$. Then
$$ \xi \cdot 2x \in \Z \tag3.2$$ for all $\xi \in \Lambda$.
\endproclaim

In the next Section we shall show how these facts can be used to
show that the only convex bodies which admit spectra are
quadrilaterals and hexagons.

\subhead Proof of Proposition 3.1
\endsubhead We may rescale so that $x = e_1$, the coordinate
direction $(1, 0 \dots, 0)$, and $I$ is the interval from $(e_1 -
e_2)/2$ to $(e_1 + e_2)/2$. Thus, we must show that
$$ \Lambda \subset \Z \times \R. \tag 3.3$$

The set $\Omega$ thus contains the unit square $Q :=
[-1/2,1/2]^2$.  Since we are assuming $\Omega$ is not a
quadrilateral, we therefore have $|\Omega| > 1$. In particular,
$\Lambda$ has asymptotic density strictly greater than $1$, i.e
the expression $(1.6)$ is strictly greater than $1$.

A direct computation shows that

$$ \hat \chi_Q(\xi_1, \xi_2) =
\frac{ \sin(\pi \xi_1) \sin(\pi \xi_2) }{\pi^2 \xi_1 \xi_2}. \tag
3.4$$ The zero set of this is
$$ Z_Q := \{ (\xi_1, \xi_2): \xi_1 \in \Z - \{0\} \hbox{ or } \xi_2 \in \Z -
\{0\} \}. \tag3.5$$
Note that $Z_Q \subset G$, where $G$ is the Cartesian grid
$$ G := (\Z \times \R) \cup (\R \times \Z). \tag3.6$$

Heuristically, we expect the zero set $Z_\Omega$ of $\hat
\chi_{\Omega}$ to approximate $Z_Q$ in the region $|\xi_1| \gg
|\xi_2|$. The following result shows that this indeed the case.
\proclaim{Lemma 3.3} For every $A \gg 1$ and $0 < \eps \ll 1$,
there exists an $R \gg A$ depending on $A$, $\eps$, $\Omega$, such
that $Z_\Omega \cap S_{A,R}$ lies within a $O(\sqrt{\eps})$
neighborhood of $Z_Q$, where $S_{A,R}$ is the slab
$$ S_{A,R} := \{
(\xi_1, \xi_2): |\xi_1| \geq R; |\xi_2| \leq A \}. \tag3.7$$
\endproclaim

\demo{Proof} Fix $A$, $\eps$. We may write
$$ \hat \chi_\Omega =
\hat \chi_{\Omega_-} + \hat \chi_Q + \hat \chi_{\Omega_+}
\tag3.8$$ where $\Omega_-$ is the portion of $\Omega$ below $x_2 =
-1/2$, and $\Omega_+$ is the portion above $x_2 = 1/2$.  In light
of (3.4), it thus suffices to show that
$$ |\hat
\chi_{\Omega_\pm}(\xi_1, \xi_2)| \lesssim \eps / |\xi_1| \tag
3.9$$ on $S_{A,R}$. By symmetry it suffices to do this for
$\Omega_+$.

We may write $\Omega_+$ as
$$ \Omega_+ = \{ (x,y): -1/2 \leq x
\leq 1/2; 1/2 \leq y \leq 1/2 + f(x) \} \tag3.10$$ where $f$ is a
concave function on $[-1/2,1/2]$ such that $f(\pm 1/2) = 0$.

By continuity of $f$, we can find a $0 < \delta \ll \eps$ such that
$$ f(1/2 - \delta), f(\delta - 1/2) \leq \eps. \tag3.11$$
Draw the line segment from $(1/2,1/2)$ to $(1/2-\delta,
1/2+f(1/2-\delta))$, and the line segment from $(-1/2,1/2)$ to
$(-1/2+\delta, 1/2+f(-1/2+\delta)$. This divides $\Omega_+$ into
two small convex bodies and one large convex body. The diameter of
the small convex bodies is $O(\eps)$, and so their contribution to
(3.9)is acceptable by Lemma 2.1. If $R$ is sufficiently large
depending on $A$, then $(\xi_1, \xi_2)$ will always make an angle
of $\gtrsim \delta/\eps$ with the normals of the large convex
body. By Lemma 2.2, the contribution of this large body is
therefore $O(\delta/\eps |\xi|^2)$, which is acceptable if $R$ is
sufficiently large.
\enddemo

Let $A \gg 1$ and $0 < \eps \ll 1$, and let $R$ be as in Lemma
3.3. Since $\Lambda - \Lambda \subset Z_\Omega$, then by Lemma 3.3
we see that $\Lambda \cap (\xi + S)$ lies in an $O(\sqrt{\eps})$
neighborhood of $Z_Q + \xi$ for all $\xi \in \Lambda$. Suppose
that we could find $\xi, \xi' \in \Lambda$ such that $|\xi - \xi'|
\ll A$ and
$$ \dist(\xi - \xi', G) \gg \sqrt{\eps}. \tag3.12$$

It follows that
$$ \Lambda \cap (\xi + S_{A,R}) \cap (\xi' + S_{A,R}) \tag3.13$$
lies in an $O(\sqrt{\eps})$ neighborhood of $G + \xi$ and in an
$O(\sqrt{\eps})$ neighborhood of $G + \xi'$.  Since $\Lambda$ has
separation $\gtrsim 1$, it follows that $\Lambda$ has density at
most $1 + O(1/A)$ in the set $(\xi + S_{A,R}) \cap (\xi' +
S_{A,R})$.  However, this is a contradiction for $A$ large enough
since $\Lambda$ needs to have asymptotic density
$1/|\Omega|<1/|Q|=1$.

By letting $\eps \to 0$ and $A \to \infty$ we see that
$$ \xi - \xi' \in G \tag 3.14$$
for all $\xi, \xi' \in \Lambda$. In particular, $\Lambda \subset
G$since $(0,0) \in \Lambda$.

Now suppose for contradiction that (3.3) failed. Then there exists
$(\xi_1, \xi_2) \in \Lambda$ such that $\xi_1 \not \in \Z$. Since
$\Lambda \subset G$, we thus have that $\xi_2 \in \Z$.  From
$(3.14)$ we thus see that
$$ \Lambda \subset \R \times \Z. \tag3.15$$

For each integer $k$, let $R_k$ denote the intersection of
$\Lambda$ with $\R \times \{k\}$.

Let $A \gg 1$ and $0<\eps \ll 1$, and let $R$ be as in Lemma 3.3.
If $\xi, \xi' \in R_k$ and $|\xi - \xi'| \gg R$, then by Lemma 3.3
we see that $\xi - \xi'$ lies in a $O(\sqrt{\eps})$ neighbourhood
of $\Z$.  From this and the separation of $\Lambda$ we see that
one has
$$ \# \{ (\xi_1, k) \in R_k: |\xi_1| \leq M \} \lesssim M+R \tag3.16$$
for all $k$ and $M$. Summing this for $-M<k<M$ and then letting $M
\to \infty$ we see that $\Lambda$ has asymptotic density at most
$1$, a contradiction. This proves $(3.3)$, and Proposition 3.1 is
proved.

\subhead Proof of Proposition 3.2 \endsubhead
By an affine
rescaling we may assume that $x = e_1/2$ and $n = e_1$, so that
our task is again to show $(3.3)$. We shall prove the following
analogue of Lemma 3.3. \proclaim{Lemma 3.4} For all $A \gg 1$, $0
< \eps \ll 1$ there exists an $R \gg 1$ depending on $A$, $\eps$,
$\Omega$ such that $Z_\Omega \cap B(Re_1, A)$ lies within
$O(\eps)$ of $\Z \times \R$. \endproclaim

\demo{Proof} Fix $A, \eps$.  We can write $\Omega$ as
$$ \Omega=
\{ (x,y): -1/2 \leq x \leq 1/2; -f(-x) \leq y \leq f(x) \}
\tag3.17$$ where $f(x)$ is a concave function on $[-1/2,1/2]$
which vanishes at the endpoints of this interval but is positive
on the interior.

For each $0<\delta \ll 1$, define
$$ S(\delta) := \frac{f(1/2 -
\delta) + f(-1/2 + \delta)}{\delta}. \tag3.18$$ The function
$\delta S(\delta)$ is decreasing to 0 as $\delta \to 0$. Thus we
may find a $0 <  \delta_0 \ll \eps/A$ such that $\delta_0
S(\delta_0) \lesssim \eps/A$.

Fix $\delta_0$, and let $l_+$, $l_-$ be the line segments from
$(1/2 -2\delta_0,0)$ to $(1/2-\delta_0, f(1/2-\delta_0))$ and
$(1/2-\delta_0, -f(-1/2+\delta_0))$ respectively, and let $-l_+$,
$-l_-$ be the reflections of these line segments through the
origin.

By symmetry we have
$$ \hat \chi_\Omega = 2 \Re (
\widehat{\chi}_{\Omega_+} + \widehat{\chi}_{\Gamma_0} )\tag 3.19$$
where $\Omega_+$ is the portion of $\Omega$ above $l_+$, $-l_-$,
and the $e_1$ axis, and $\Gamma_0$ is the small portion of
$\Omega$ between $l_+$ and $l_-$.

Since we are assuming $\Omega$ to have normal $e_1$ at $e_1/2$, we
see that $S(\delta) \to \infty$ as $\delta \to 0$. Thus we may
find a $0 < \delta \ll \delta_0$ such that
$$ S(\delta) \gg 1 + \frac{1}{\eps} S(\delta_0). \tag3.20$$

Fix this $\delta$. By Corollary 2.4 we may find an $R \sim
1/\delta$ such that
$$ |\hat \chi_{\Gamma_0}(R e_1)| \gtrsim (f(1/2
- \delta) + f(-1/2 + \delta)) \delta = \delta^2 S(\delta). \tag
3.21$$

Fix this $R$. Let $m_+$, $m_-$ be the line segments from $(1/2 -
2\delta, 0)$ to $(1/2 - \delta, f(1/2 - \delta))$ and $(1/2 -
\delta_0, -f(-1/2 + \delta))$ respectively. We can partition
$$ \hat \chi_{\Gamma_0} = \hat \chi_{\Gamma_+} + \hat \chi_{\Gamma_-}
+ \hat \chi_\Gamma \tag3.22$$ where $\Gamma_+$ is the portion of
$\Gamma$ above $m_+$ and the $e_1$ axis, $\Gamma_-$ is the portion
below $m_-$ and the $e_1$ axis, and $\Gamma$ is the portion
between $m_+$ and $m_-$.

The convex body $\Gamma - e_1/2$ is contained inside a ball of
radius $O(S(\delta) \delta)$, hence by (0.2) we have

$$ |\nabla
\hat \chi_{\Gamma - e_1/2}(\xi)| \lesssim (\delta S(\delta))^2 / R
\lesssim (\delta_0 S(\delta_0)) \delta^2 S(\delta) \lesssim
\frac{\eps}{A} \delta^2 S(\delta) \tag3.23$$ for $\xi \in B(Re_1,
A)$.

If $\xi \in B(Re_1, A)$, then $\xi$ makes an angle of
$$O(A/R)=O(A\delta) \ll O(\delta/(\delta_0 S(\delta_0)))=
O(\delta/(\delta S(\delta))) = O(1/S(\delta)) \tag3.24$$ with the
$e_1$ axis, and hence makes an angle of $\gtrsim 1/S(\delta)$ with
the convex bodies $\Gamma_+ - e_1/2$, $\Gamma_- - e_1/2$. Since
these bodies are in a ball of radius $O(S(\delta_0)
\delta_0)=O(\eps/A)$, we see from Lemma 2.2 that
$$ |\nabla \hat
\chi_{\Gamma_\pm - e_1/2}(\xi)| \lesssim \frac{\eps}{A} S(\delta)
/ R^2 \sim \frac{\eps}{A} \delta^2 S(\delta). \tag3.25$$

Summing, we obtain
$$ |\nabla (e^{\pi i \xi_1} \hat
\chi_{\Gamma_0}(\xi))| \lesssim \frac{\eps}{A} \delta^2 S(\delta).
\tag3.26$$

Integrating this and $(3.21)$ we get
$$ \hat \chi_\Gamma(\xi) =
\hat \chi_{\Gamma_0}(R e_1)( e^{\pi i (R-\xi_1)}+O(\eps) ).\tag
3.27$$

If $\xi \in B(Re_1, A)$, then $\xi$ makes an angle of $\gtrsim
1/S(\delta_0)$ with every normal of $\Omega_+$. From Lemma 2.2 we
get
$$ |\hat \chi_{\Omega_+}(\xi)| \lesssim  S(\delta_0) / R^2 \sim
S(\delta_0) \delta^2 \ll \eps \delta^2 S(\delta) \tag3.28$$ on
$B(Re_1, A)$.  From this, $(3.20)$, $(3.22)$, and $(3.23)$ we
obtain
$$ \hat \chi_\Omega(\xi) = 2 \Re( \hat \chi_\Gamma(R e_1)(
e^{\pi i(R - \xi_1)}+O(\eps) )) \tag3.29$$ on $B(Re_1, A)$, and
the Lemma follows.
\enddemo

Let $A \gg 1$, $0 < \eps \ll 1$, and let $R$ be as in Lemma 3.4.
If $\xi \in \Lambda$ are such that $|\xi| \ll A$, then from Lemma
3.4 we see that
$$ Z_\Omega \cap B(Re_1 + \xi, A) \cap B(Re_1, A)
\tag3.30$$ lies within $O(\eps)$ of $(\Z \times \R)$, and within
$O(\eps)$ of $(\Z \times \R) + \xi$. Since $Z_\Omega$ has
asymptotic density $1/|\Omega|$, it has a non-empty intersection
with $B(Re_1, A) \cap B(Re_1 + \xi, A)$, and thus $\xi$ must lie
within $O(\eps)$ neighbourhood of $\Z \times \R$. Taking $\eps \to
0$ and then $A \to \infty$ we obtain $(3.3)$, and Proposition 3.2
is proved.

\head Conclusion of the argument \endhead \vskip.125in

We now use Proposition 3.1 and 3.2 to show that the only convex
symmetric bodies with spectra are the quadrilaterals and hexagons.
We may assume of course that $\Omega$ is not a quadrilateral or a
hexagon.

Suppose that there are two points $x$, $x'$ in $\partial \Omega$
for which either Proposition 3.1 or Proposition 3.2 applies. From
elementary geometry we thus see that $\Lambda$ must live in a
lattice of density $|2x \wedge 2x'|$. It follows that

$$ 4|x \wedge
x'| \geq |\Omega| \tag4.1$$ for all such $x, x'$.  Since $|x| \sim
1$ on $\partial \Omega$, this implies that there are only a finite
number of $x$ for which Proposition 3.1 and Proposition 3.2
applies. Since almost every point in $\partial \Omega$ has a unit
normal, the only possibility left is that $\Omega$ is a polygon.

Label the vertices of $\Omega$ cyclically by $x_1, \ldots,
x_{2n}$. Since $\Omega$ is not a quadrilateral or a hexagon, we
have $n \geq 4$. By symmetry we have $x_{n+i} = -x_i$ for all $i$
(here we use the convention that $x_{2n+i}=x_i$).

From Proposition 3.1 we have
$$ \xi \cdot (x_i - x_{n+i-1}) \in \Z
\tag4.2$$ for all $\xi \in \Lambda$. First suppose that $n$ is
even. Then $n-1$ is coprime to $2n$, and by repeated application
of $(4.2)$ we see that
$$ \xi \cdot (x_i - x_j) \in \Z \tag4.3$$
for all $i, j$.  Arguing as in the derivation of $(4.1)$ we thus
see that
$$ |(x_i - x_j) \cdot (x_i - x_k)| \geq |\Omega| \tag4.4$$
for all $i, j, k$.  In other words, the triangle with vertices
$x_i$, $x_j$, $x_k$ has area at least $|\Omega|/2$ for all $i, j,
k$.  But $\Omega$ can be decomposed into $2n-2$ such triangles, a
contradiction since $n \geq 4$.

Now suppose that $n$ is odd, so that $n \geq 5$. Then $n-1$ and
$2n$ have the common factor of 2. Arguing as before we see that
$(4.3)$ holds for all $i,j,k$ of the same parity. But $\Omega$
contains the three disjoint triangles with vertices $(x_1, x_3,
x_5)$, $(x_1, x_5, x_7)$, and $(x_1, x_7, x_9)$ respectively, and
we have a contradiction.

\newpage
\vskip.25in
\head References \endhead

\ref \key Fug74 \by B. Fuglede \paper Commuting self-adjoint partial
differential operators and a group theoretic problem \jour J. Funct.
Anal. \yr 1974 \vol 16 \pages 101-121 \endref

\ref \key Herz64 \by C. S. Herz \paper Fourier transforms related to
convex sets \jour Ann. of Math. \vol 75 \yr 1962 \pages 81-92
\endref

\ref \key IKP99 \by A. Iosevich, N. Katz, and S. Pedersen
\paper Fourier basis and the Erd\H os distance problem
\jour Math. Research Letter \yr 1999 \vol 6 \pages \endref

\ref \key IKT01 \by A. Iosevich, N. Katz, and T. Tao \paper Convex
bodies with a point of curvature do not have Fourier bases \jour
Amer. J. Math. \vol 123 \yr 2001 \pages 115-120 \endref

\ref \key IosPed99 \by A. Iosevich and S. Pedersen
\paper How wide are the spectral gaps? \jour Pacific J. Math.
\vol 192 \yr 1998 \pages 307-314 \endref

\ref \key Kol00 \by M. Kolountzakis \paper Non-symmetric convex
domains have no basis of exponentials \jour Illinois J. Math. \vol
44 \yr 2000 \pages 542--550 \endref

\ref \key Lab01 \by I. Laba \paper "Fuglede's conjecture for a
union of two intervals (preprint) \yr 2001 \endref

\ref \key Lan67 \by H. Landau \paper Necessary density conditions for
sampling and interpolation of certain entire functions \jour Acta Math.
\vol 117 \pages 37-52 \yr 1967 \endref

\enddocument

\enddocument